%% file: weillinois.tex
\documentstyle[12pt]{article}

\input def.tex

 \title{{\sc Orbits of conditional expectations}}
 \author{M. Argerami and D. Stojanoff}
 \date{}

\begin{document}

 \maketitle
 
 %\renewcommand{\theequation}{\arabic{section}.\arabic{equation}}
 %\renewcommand{\theenumi}{\roman{\enumi}}

 %\begin{abstract}{
 %{\sc Abstract}
 
%\subsection*{Abstract.}
\begin{abstract}
Let $N \subseteq M$ be von Neumann algebras and        
$E:M\to N$ a faithful normal 
conditional expectation. In this work it is shown that the 
similarity orbit ${\cal S}(E)$ of $E$ by the natural action of the 
invertible group of $G_M$ of $M$ has a natural complex analytic structure
and the map given by this action: $G_M\to {\cal S}(E)$  is a smooth principal bundle. 
It is also shown that if $N$ is finite then ${\cal S}(E)$ admits a Reductive 
Structure. 
These results were known previously under the conditions of finite index and
$N'\cap M \subseteq N$, which are removed in this work. 
Conversely, if the orbit ${\cal S}(E)$ has an Homogeneous Reductive 
Structure for every expectation defined on $M$, then $M$ is finite. 
For every algebra $M$ and every expectation $E$, 
a covering space of the unitary 
orbit ${\cal U}(E)$ is constructed in terms of the connected 
component of $1$ in the normalizer  of $E$. 
Moreover, this covering  space is the universal covering in any of the 
following cases: 1) $\m$ is a finite factor and $Ind(E) < \infty $;
2) $M$ is properly infinite and $E$ is any expectation;
3) $E$ is the conditional expectation onto the centralizer of a state.
Therefore, in those cases, the fundamental group of 
${\cal U}(E)$ can be characterized as the Weyl group of $E$. 
\end{abstract}

 \section{Introduction.}

 Let $\m$ be a von Neumann algebra with group of invertible elements $\gm$ and
 unitary group $\um$. Denote by $\em$ the space of faithful
 normal conditional expectations defined on $\m$ and  ${\cal B}(\m)$ the
 algebra of bounded linear operators on $\m$. 
 Consider the action 
 \[
 L:\gm\times {\cal B}(\m)\ra {\cal B}(\m) 
 \] 
 given by
 \[ 
 L_g(T)=gT(g^{-1}\cdot g)g^{-1}, \quad g \in \gm , \ T \in {\cal B}(\m). 
 \]
 Let $\E \in \em$ be a conditional expectation. Define 
 the unitary orbit of $\E$ 
 \begin{equation}\label{defouni}
 \uE=\{L_u(\E):u\in\um\} \inc \em .
 \end{equation}
 considered with the quotient topology induced by the norm topology of $\um$.
 So we have a natural fibration 
 \begin{equation}\label{pisube}
 \Pi_\E : \um \to \uE \quad \hbox { given by } \quad 
 \Pi_\E (u) = L_u(\E), \ u\in \um 
 \end{equation}
 The objetive of this work is the study of the homotopy groups and the 
 differential geometry of the orbit $\uE$ or, more precissely, of the 
 fibration $\Pi_\E$. Several results in this sense appear in 
 the works \cita {as}, \cita{alrs} and \cita{pucho}, mainly under two 
 very restrictive hypothesis: finite index condition for $\E$ and the 
 condition $\E (\m)'\cap\m\inc\E (\m )$ for the inclusion of the 
algebras. Also in \cita {av2} a complete study of this problems is made in  
 the case when $\E$ is a state.

 In order to study the homotopy type of these orbis we construct
 a covering space over each orbit $\uE$ whose 
 group of covering trnsformations is the so called Weyl group of 
 the expectation $E$. To describe this structure we need 
 the following definitions:
 
 \bigskip
 At the level of the unitary group $\um $ of $\m$, the isotropy group of
 the action $L$, i.e. $\Pi_\E \inv (\E)$ is
 a very well known group usually called the {\it normalizer} of $\E$. It  
 has already been studied, between other authors, 
 by A. Connes (\cita{connes}) and 
 Kosaki (\cita{kosaki2}) in relation with crossed 
 product inclusions of algebras (see also \cita {ars}). 
 We shall denote the normalizer of $\E$ by 
 \begin{equation}\label{nE}
 \nE =\ \{\ u\in\um \ : \ \E(uxu^*)=u\E(x)u^*, x\in\m \ \}.
 \end{equation}
 
 Let $\n = \E(\m)$. Then $\n $ is a von Neumann 
 subalgebra of $\m$. Consider also the \avn 
 \begin{equation}\label{me}
 \me = \{ x \in N'\cap M \ : \ E(xm) = 
 E(mx) \ \hbox{ for all } \ m \in M \},
 \end{equation}
 usually denoted as the centralizer of $E$ (see \cita{combes} or \cita{h}). 
 In \cita {ars} it is shown that the connected 
 component of $1$ in $\nE$ is the group 
 \begin{equation}\label{he}
 \he =\un \cdot \ume =  \{ vw: v \in \un \; \hbox{ and } 
 \; w \in \ume \},
 \end{equation}
 which is a close, open and invariant subgroup of $\nE$. Then 
 the set of connected components of $\nE$ is a discrete group 
 called the Weyl group of $\E$:
 \begin{equation}\label{we}
 \we = \pi_0 (\nE) \simeq \nE/\he ,
 \end{equation}
 which has several characterizations in very different contexts 
 (see \cita{kosaki2}, \cita{ars} and \cita{ars2}).

\bigskip
 We show that, for any \avn \ $\m$ and any $\E\in \em$, 
 the space $\ee = \um /\he$ and its natural 
 projection onto $\uE \simeq \um/\nE$ 
 (see diagram (\ref{diagrama1})) defines a covering map whose group of
 covering transformations can be identified with the 
 Weyl group $ \we$ 
(Theorem \ref{covering}).

In all examples we know, 
$\ee$ is actually the $universal$ covering for $\uE$ and therefore the 
fundamental group $\pi_1 (\uE , \E) $ coincides with the Weyl group $\we$.
We conjeture that this is true for all \avns \ $\m$ and all conditional
expectations $\E\in \em$. In Theorem \ref{uni} we show that if any of the 
following conditions hold:
\ben
\item $\m$ is properly infinite, 
\item $\m$ is finite, $\dim \zm <\infty $ and $\E$ has finite index,
\item $\E= \E_\varphi \in {\cal E} (\m , \m _\varphi )$ is the canonical 
expectation associated to a faithful normal state $\varphi$ of $\m$, i.e. 
$\m_\varphi $ is the centralizer of $\varphi$ and $\varphi \circ \E =
\varphi$, 
\een

\noindent 
then $\ee$ is simply connected, so that it is the universal covering for the 
orbit $\uE$. Consequently,
$$\pi_1(\uE , \E )\simeq\we.$$

\bigskip
In order to study the differential geometry of the orbit of an expectation
$\E$ we consider the whole similarity orbit 
\begin{equation}\label{defoinvertible}
\o=\{L_g(\E):g\in\gm\} \inc  {\cal B}(\m) 
\end{equation}
and the fibration (with the same name as its restriction to $\um$), 
\begin{equation}\label{gpisube}
\Pi_\E : \gm \to \oe \quad \hbox { given by } \quad 
\Pi_\E (g) = L_g(\E), \ g\in \gm .
\end{equation}
Note that $L_g(\E)$ is not necessarily a conditional expectation for all 
$g \in \gm$. Nevertheless we prefer to use 
this setting, since the group $\gm$ 
is a complex analytic Banach Lie group and then the orbit $\o$ can be 
given a complex analytic manifold structure. In any case, all the 
geometrical results obtained for $\o$ are also true for the unitary 
orbit $\uE$ by just replacing ``complex analytic'' by ``real analytic''.

\bigskip
In order to study the differential geometry of similarity orbits we need 
to generalize several results of \cita {ars} mentioned before 
to the invertible groups setting. This task is made
in section 3, where the connected component of the isotropy group
$\ie$ of the action $L$ at $\E$ is characterized 
(see Proposition \ref{conexa}) 
and the new Weyl group which naturally appears is shown to be the 
same group as the ``old'' one (see Theorem \ref{W(E)}).

\bigskip

In section 4 we first show that $\o$ can be always be given
a unique complex analytic differential structure such that the map 
$\Pi_\E$ of equation (\ref{gpisube}) becomes a submersion 
(Theorem \ref{estructura}). The key tool is the construction, in the
style of \cita {combes}, of a
conditional expectation $\F\in \em$ onto the centralizer $\me$
which commutes with $\E$. This allows us to get a complement in $\m$ 
for the subspace $\n + \me$, which can be naturally identified with
the tangent space of $\o$ at $\E$. 

Then we show that if $\n = \E(\m) $ is a finite \avn , $\o$ has a 
unique structure of Homogeneous Reductive Space (HRS) 
(see Definition \ref{HRS} and Proposition \ref{HRS1}). 
This family of HRS's are of geometrical interest. Indeed, 
maybe the most general families of examples of infinite dimensional 
HRS's modeled in operator algebras are studied in \cita{acs1}
and \cita{alrs}. But all these examples can be 
represented as the quotient of the unitary (or invertible) group of an 
algerba $\m$ by the unitary (or invertible) group of some subalgebra. 
This situation does not happen in the case of the orbit
of a conditional expectation. Indeed, the isotropy group $\ie$ can be   
big enough even to generate the whole algebra $\m$, while $\o \simeq
\gm / \ie$.  Moreover, also the map $\Pi_0 : \gm \to \gm / \ze $ ($=\eee$), 
where $\ze$ is the connected component of $1$ in $\ie$, defines a HRS if
$\n$ is finite. Actually (see Teorem \ref{estructura})
this is the way to show the HRS structure
of $\o$, since $\eee$ is a covering space for $\o$ and therefore they are
locally homeomorphic (and diffeomorphic). But neither $\eee$ can be 
represented as a quotient as before (if $\n \not\inc  \me $ and 
$\me \not\inc \n$), since, by Proposition \ref{mdecomp}, 
$\ze = \gme \cdot \gn$ and this is not the invertible group of any subalgebra 
of $\m$.

In the end of section 4 we show 
that the existence of HRS structures for any expectation $\E\in \em$ 
forces $\m$ to be a finite \avn \ (Theorem \ref{Principal}).

\section{The universal covering of $\uE$.}
Let $\n\subseteq\m$ be von Neumann algebras. From now on we shall 
denote by $\emn$ the space of 
faithful normal conditional expectations $\E: \m \to \n$. 
Let $\E \in \emn $. We recall the definitions of
the sets $\uE $, $\nE$, $\me$ and $\he$ (see equations 
(\ref {defouni}), (\ref{nE}), (\ref{me}) and (\ref{he}), 
respectively) associated with $\E$. Consider the space $\ee = \um / \he$, 
with the quotient topology of the norm topology of $\um$ and denote by 
$\Pi_0$ the projection form $\um$ onto $\ee$. 
The situation we shall study is the following: we have 
a commutative diagram
\begin{equation}\label{diagrama1}
\begin{array}{rcl}
\um&\hrarr{\Pi_0}{}&\ee=\um/\he  \\
&\ddrarr{\Pi_\E}{}&\vdarr{}{\Phi}\\
&& \uE \simeq  \um/\nE
\end{array}
\end{equation}
where the map $\Phi$ is defined by 
$\Phi (u \he ) = \Pi_\E (u) \sim u \nE $, $ u\in \um$.
In \cita{as} it was shown that when $\n'\cap \m\inc\n$, and the 
Jones index of $\E$ is finite, calling $e$ the Jones projection of 
$\E$, then its $\m$-unitary orbit  
$$
{\cal U}_M(e) = \{ueu^* : u \in \um \} \simeq \um / \un ,
$$
is a covering space for $\uE$. Note that, under the mentioned assumptions, 
${\cal U}_M(e) \simeq \ee$, since both spaces can be  indentified
with $\um / \un $ and the fact, showed in \cita {as}, that the qoutient 
topology and the norm topology coincide on ${\cal U}_M(e) $ .

In what follows we will show, removing both hypothesis appearing 
in \cita {as}, that the map $\Phi$ is always a covering map, and 
that $\ee = \um/\he$
is a covering space for $\uE$, with group of covering transformations $\we$.
Moreover, in several cases (see \ref{uni}) $\Phi$ is the universal covering 
of $\uE$ and, in particular, $\pi_1(\uE)\simeq\we.$

\bigskip
Note that the Weyl group $\we= \nE /\he $, being included in $\ee$, 
has a natural action on $\ee$ given by right
multiplication. This action is well defined 
because $\he$ is a normal subgroup of $\nE$.

\begin{pro}\label{fibra}
Let $\n\subseteq\m$ be von Neumann algebras and $\E \in \emn $
a faithful normal conditional expectation. 
Then, with the notations of diagram (\ref{diagrama1}), 
\ben
\item The map $\Phi$ is continous.
\item For any $u \in \um$, the fibre by $\Phi$ of $L_u\E \in \uE$ is 
precisely the orbit of  $\Pi_0(u) \in \ee$ by the 
action of $\we$.
\item The unitary orbit $\uE$ is homeomorphic 
with $\ee/\we$ (i.e. the space of orbits by the action of $\we$ in $\ee$), 
both considered with the quotient topology.
\een
\end{pro}
\dem \ Items 1 and 2 follow immediately by looking at 
the commutative diagram (\ref{diagrama1}) and 
using the fact that $\Phi \inv (\E) = \we$. 
Let $\rho:\ee \to \ee/\we$ 
be the canonical projection. 
Consider the map $\bar{\Phi}:\ee/\we \to \uE $ 
given by $\bar{\Phi}(\rho(h))=\Phi(h)$, 
$h\in \ee$. Then $\bar{\Phi}$ is the desired homeomorphism.
Indeed, it is clear that $\bar{\Phi}$ is well defined and  is 
bijective. $\bar{\Phi}$ is also continuous, since 
$\Phi=\bar{\Phi} \circ \rho $ is continuous. On the other hand, 
let $U$ be an open set in $\ee/\we$. By considering the full commutative 
diagram
\begin{equation}
\begin{array}{rcccl}
\um&\hrarr{\Pi_0}{}&\ee=\um/\he&\hrarr{\rho}{}&\ee/\we\\
&\ddrarr{\Pi_\E}{}&\vdarr{\Phi}{}&\ddlarr{\bar{\Phi}}{}&\\
&&\uE=\um/\nE&&
\end{array}
\end{equation}
and the fact that $\Pi_\E$ is an open map, 
%(since it  is a submersion, see Theorem \ref{estructura} below), 
it is clear that $\bar{\Phi}(U)$ is open in $\uE$.\qed

\begin{rem}\label{teogreenberg}\rm
In order to show that the map $\Phi$ defined in diagram (\ref{diagrama1}) 
is a covering map we shall use the following well known result of algebraic
topology (see, for instance, Chapter 1 of \cita{greenberg}):

Let $X$ be a locally pathwise connected and connected topological 
space, and $G$ a 
groups of homeomorphisms of $X$ that operates properly discontinously (i.e.
for each $x\in X$ there exists an open set $V_x$ such that 
$V_x\cap g(V_x) =\empty$ for
every $g\in G$, $g\ne e$). Consider the map $p:X\ra X/G$. 
Then $X$ is a covering space,
for $X/G$ with covering map $p$ and group of covering transformations $G$, 
and $p_*(\pi_1(X, x_0))$ is a normal subgroup of $\pi_1(X/G, p(x_0))$.
\end{rem}

\begin{teo} \label{covering}
Let $\n\subseteq\m$ be von Neumann algebras and $\E \in \emn $
a faithful normal conditional expectation. 
Then, with the notations of diagram (\ref{diagrama1}), the space $\ee$ 
is a covering space for $\uE$, 
with covering map $\Phi$ and group of covering transformations $\we$.
\end{teo}
\dem By the previous remark, it suffices to show that $\we$ 
operates properly discontinously on $\ee$. Fix $u \in \um$ and consider 
the open set 
$$
W_u=\{ w \in \um : \|w-u\|< 1/2 \}.
$$
For each element $k \in \we$, we  choose some 
$u_k \in \nE$ such that $\Pi_0(u_k) = u_k \he = k$. 
Since the map $\Pi_0$
%, being a submersion (see \ref{lie3} and \ref{mdecomp} below), 
is open, we can consider the open set 
$$
V_u = \Pi_0 (W_u) \inc \ee.
$$
Note that, for $k \in \we$, $V_u k = \Pi_0 (W_u u_k)$. 
In order to prove that the action of $\we$ in $\ee$ is
properly discontinous, we just have to see that 
$V_u \cap V_u k = \phi$  for every $1\ne k \in \we$. 
Suppose that this is not true. Then, for some $1\ne k \in \we$, 
there exist $w_1, w_2 \in W_u$ and 
$ z \in \he$ such that $w_1 u_k = w_2 z$. Then 
$$
 w_1^* w_2 = u_k z^*  \in \nE \setminus \he
$$
But, since $w_1, w_2 \in W_u$,  
\[
\|w_1^* w_2 -1 \| =  \|w_2 - w_1 \|  < 1.
\]
This implies that $w_1^* w_2 \in \he$ (see \cita{ars} or the proof of  
Proposition \ref{conexa} below), a contradiction.
\qed

\begin{cor}\rm\label{isopi}
The group $\Psi_*(\pi_1(\ee))$ is a normal subgroup of $\pi_1(\uE)$ 
and we have the following isomorphism:
\[\pi_1(\uE)/ \Psi_*(\pi_1(\ee)) \simeq  \we \]
\end{cor}
\dem By Proposition \ref{fibra}, we know that the fibre $\Psi^{-1}(\E)$ 
equals $\we$, and the assertion follows by using the homotopy exact sequence 
induced by the covering map $\Psi $. \qed

\begin{rem}\rm \label{ava}
Let $\varphi $ be a faithful normal state of the \avn \ $\m$. In
\cita {av2} Andruchow and Varela  show that the unitary orbit of $\varphi$: 
$$
{\cal U}(\varphi) = \{ \varphi ( u^* \cdot u ) : u \in \um \}
$$
is simply connected.  Therefore the unitary group of the 
centralizer $\m _\varphi $ of $\varphi $ coincides with the  
normalizer ${\cal N}_\varphi $ of $\varphi$, considered as a 
conditional expectation. Then the covering space 
$$
{\cal X }(\varphi ) =\um /  {\cal U}_{\m _\varphi} 
=\um /  {\cal N}_{\varphi} \simeq {\cal U}(\varphi) 
$$ 
and ${\cal U}(\varphi)$ is its own universal covering.

Moreover, if $\E_\varphi \in {\cal E}(\m , \m _\varphi )$ is the canonical 
expectation such that $\varphi \circ \E _\varphi = \varphi $, then 
${\cal U}(\varphi ) \simeq {\cal X} (\E _\varphi )$ and so it is the
universal covering for ${\cal U}(\E_\varphi )$. Indeed, since
$${\cal X} (\E _\varphi ) = \um / {\cal U}_{\m _\varphi} 
{\cal U}_{\m _{\E_\varphi}} \quad \mbox{ and }
 \quad {\cal U}(\varphi) \simeq 
\um /  {\cal U}_{\m _\varphi }, 
$$
it suffices to show that $\m _{\E_\varphi} \inc \m _{\varphi }$. 
But this is apparent by the definition of $\m _{\E_\varphi}$  (see equation
\ref{me}) and the fact that $\varphi \circ \E _\varphi = \varphi $.
This gives a large class of conditional expectations for which the
covering space $\ee$ is the universal covering. We shall extend this 
class in the following Theorem.
\end{rem}

\begin{teo}\label{uni}
Let $\m$ be a \avn , $\E\in\em$, and suppose that 
any of the following conditions hold:
\ben
\item $\m$ is properly infinite. 
\item $\m$ is a II$_1$ factor and $\E$ has finite index.
\item $\E= \E_\varphi \in {\cal E}(\m , \m _\varphi )$ is the canonical 
expectation associated to a faithful normal 
state $\varphi$ of $\m$  as in Remark \ref{ava}.
\een

\noindent 
Then $\ee$ is simply connected, so that it is the universal covering for the 
orbit $\uE$. Consequently,
$$\pi_1(\uE , \E )\simeq\we.$$
\end{teo}
\dem
Consider the fibre bundle 
\begin{equation}\label{fibre}
\Pi_0:\um \to \um /\he = \ee. 
\end{equation}

\noindent
Recall that a fibre bundle gives rise to an exact 
sequence of homotopy groups. 
In our case, the bundle $\Pi_0$ yields the exact sequence
\begin{equation}\label{larga}
\ldots\to \pi_2 (\ee) \to \pi_1(\he) \stackrel{i_*}{\to} \pi_1 (\um) \to 
\pi_1 (\ee) \to \pi_0 (\he) = 0 ,
\end{equation}
where $1$ is taken as base point for the homotopy groups of 
the unitary groups and $[1]_{\ee}= \he$ is the base point for $\ee$. Here 
$i_*$ denotes the homomorphism induced by the inclusion 
$i : \he \hookrightarrow \um$.
We can then use results by Handelmann \cita{han} and Schr\"oeder 
\cita{sch} on computing the homotopy group of
the unitary group of a von Neumann algebra.

\smallskip
\noindent
Case (1.) It follows by appealing to the homotopy exact sequence (\ref {larga}), 
and  the fact \cita{han} that $\um $ has trivial $\pi_1$ group if $\m$ is 
properly infinite.

\smallskip
\noindent
Case (2.): Since $\m$ is a II$_1$ factor and $\ind < \infty $ it is known 
(see \cita{Po}) that $\n =\E(\m)$ is 
also of type II$_1$ and $\dim \zn < \infty$. 
Let us recall the following results (see \cita{av2}, 
\cita{han} and \cita{sch}):
\ben
\item If $\m$ is a \avn \ of type II$_1$, then 
$\pi_1(\um)$ is isomorphic to (the additive group) 
$\zm_{sa}$ of selfadjoint elements in $\zm$. 
\item Let $j : \un \to \um $ be the inclusion map. Then 
the image of  the homomorphism $j_* :\pi_1 (\un) \to \pi_1 (\um)
\simeq \zm_{sa}$ is equal to the additive group generated by
the set $\{tr(p) : p \hbox{ projection in } \n \}$, where $tr$ is 
the center valued trace of $\m$.
\een
In our case  $\pi_1(\um) \simeq \zR$. 
Let $k : \un \to \he$ be the inclusion map.
Clearly $i_* \circ k_* = j_*$, where $i_*$ is the map of equation 
(\ref{larga}). Then $j_* (\pi_1(\un )) \inc i_*(\pi_1(\he )) $. 
Let $p\in \zn$ be a minimal 
projection. Then $p\n p$ is a II$_1$ factor and the trace of projections
only in $p\n p$ generate the additive group $\zR$. Then $i_*$ is surjective
and $\pi_1 (\ee)$ must be trivial by the homotopy exact sequence 
(\ref{larga}). 
\smallskip
\noindent
Case (3.) follows from \cita{av2} and Remark \ref{ava}  \qed
\medskip
\begin{rem}\rm
With the same techniques used in the proof of the previous Theorem, 
it can be also shown that $\ee$ is simply connected if $\ind <\infty$, 
$\m$ is finite and $\dim \zm <\infty$. 
\end{rem}

\begin {exa}\rm
Let $\m$ be a \avn \ and $p\in \m $  be a projection. 
Then $p$ determines the conditional expectation 
$\E_p : \m \to \n = \{p\}'\cap \m $ given by
$$
\E_p (x) = pxp +(1-p)x(1-p) , \quad x \in \m .
$$
Denote by $\up = \{upu^* : u \in \um\}$ the unitary orbit of $p$, which
is a connected component of the Grasmmanians of $\m$. Then 
$$
\up \simeq \eep 
$$ 
in the sense that both spaces are homeomorphic to $\um /\un $, since 
$\n ' \cap \m \inc \n$ and so $\hep = \un$. Note that on $\up$ 
we consider the norm topology as a subset of $\m$ (see \cita{cpr} or \cita{PR}).
Using Theorem \ref{uni} it is not difficult to show that the 
Grasmmanian $\up$ is always simply
connected. Indeed, $\pi_1 (\up )$ splits in the finite and the
properly infinite parts of $\m$ and  items 1 and 3 of \ref{uni} can by 
applied (see also \cita{acs3}). 

The Weyl group of $\E_p$ is trival 
if $1-p \not \in \up $ and it has two elements if $1-p \in \up$, since 
in this case, if $u \in \um $ verifies that $upu^* = 1-p$, then 
$L_u \E_p = \E_{1-p} = \E_p$ and $\nep = \un \cup u \cdot \un $.

\medskip 
A similar study can be done for sistems of projections, i.e the $n$-tuples 
$P= (p_1, \dots , p_n) $ of pairwise orthogonal projections such that
$\sum p_i = 1$ (see \cita{cpr2}). Here also (by Thm. \ref{uni}) 
the joint unitary orbit 
of $P$ is simply connected and it is homeomorphic to the space
${\cal X}(E_P)$ associated to the conditional expectation
$E_P(x) = \sum p_i x p_i $, $x \in \m$. The Weyl group is a subgroup of
the permutation group $S_n$, determined by those projections in $P$ 
such that are  equivalent (and therefore unitary equivalent) in $\m$.
\end{exa}

\section{The Weyl group, invertible case.}

\

In \cita{ars}, the Weyl group is defined in terms of the unitary group of
the von Neumann algebra $\m$. Our aim in this section is to extend 
the previous development to the case
when the action over the conditional expectations is given not anymore by
the unitaries but by the invertible elements.

Let us recall some definitions. 
Let $\m$ be a von Neumann algebra. We consider the action 
$ L:\gm\times {\cal B}(\m)\ra {\cal B}(\m) $ given by
$ L_g(T)=gT(g^{-1}\cdot g)g^{-1}$, $g \in \gm$,  $T \in {\cal B}(\m)$. 
Let $\E \in \em$ be a conditional expectation. Then,
as we have already mentioned,  $L_g(\E)$ 
is not necessarily a conditional expectation for all $g \in \gm$, but we  
still consider the orbit of the expectation
$$
\o=\{L_g(\E):g\in\gm\}
\quad \mbox{ and the fibration } \quad 
\Pi_\E : \gm \to \oe 
$$
given by $\Pi_\E (g) = L_g(\E), \ g\in \gm$. 
The role played in the unitary case by the normalizer 
$\nE=\{u\in\um: L_u(\E)=\E\}$ 
as the isotropy group of the action is now played by
\begin{equation}\label{ie}
\ie=\{g\in\gm: L_g(\E)=\E\}, 
\end{equation}
Let $\n = \E (\m) \inc \m$. Then $\n$ is a  \avn . Recall that the
centralizer of $\E$ is the \avn \ 
$\me = \{ x \in N'\cap M \ : \ E(xm) = E(mx)$ for all  $m \in M \}$.
Define the group 
\begin{equation}\label{ze}
\ze=\gme\cdot\gn \inc \ie 
\end{equation}

\begin{pro}\label{lema1}
Let $\m$ be a \avn , $\E \in \em$ and consider the groups $\ie$ 
and $\ze$ defined in equations (\ref{ie}) and (\ref {ze}). Then
\begin{enumerate}
\item\label{lema1_1} If $g\in\ie$, then $g\E(g^{-1})=\E(g^{-1})g\in \me$.
\item\label{lema1_2} If $g\in\ie$ and $\E(g^{-1})$ is invertible, 
then $g\in\ze$.
\item\label{lema1_3} $\ie\cap\m^+=\ze^+$.
\end{enumerate}
\end{pro}
\dem If $g\in\ie$, then $g\E(g^{-1})=g\E(g^{-1}g^{-1}g)=\E(g^{-1})g$. 
Let $\n = \E (\m)$ and  $b\in\n$. Then 
\[
g \E(g^{-1})b=g\E(g^{-1}b)=g\E(g^{-1}bg^{-1}g)=\E(bg^{-1})g=b\E(g^{-1})g,
\]
so that $g\E(g^{-1})\in\n'\cap \m $. If $x\in\m$, since $g\n g\inv =\n$,
\[
\E(g\E(g^{-1})x)=\E(\E(g^{-1})gx)=\E(g^{-1})\E(gx)=
\E(g^{-1})g\E(xg)g^{-1}=\E(xg\E(g^{-1})),
\]
thus proving that $g\E(g^{-1})\in\me$. 

If $\E(g^{-1}) $ is invertible, then 
$g = g\E(g^{-1}) \cdot \E(g^{-1})\inv \in \gme \gn = \ze$ by 1.
Finally, if $g\in\ie$ and $g>0$, then $g^{-1}>0$, and, as $\E$ is faithful, 
$\E(g^{-1})>0$. Then $E(g^{-1})\in\gn$, and by \ref{lema1_2}, $g\in\ze$. \qed
\medskip
\begin{lem}\label{inversa}
If $g\in\gm$ and $\|g-1\|<\eps<1$, then \[\|g^{-1}-1\|<\frac{\eps}
{1-\eps}.\] In particular, if $\eps\le 1/2$, then 
\[\|g^{-1}-1\|<2\eps.\]
\end{lem}
\dem Straightforward. \qed
\medskip
\begin{pro}\label{conexa}
Let $\m$ be a \avn , $\E\in\em$ and cosider the groups $\ze \inc \ie$ 
defined in equations (\ref{ie}) and (\ref {ze}).
The group $\ze$ is open, closed and connected in $\ie$.
Moreover, the connected component of $\ie$ at any
$u\in\ie$ is exactly $u\cdot \ze$. 
\end{pro}
\dem First of all we will show that $\ze$ is open at $1$. Let $g\in\ie$ such
that $\|g-1\|<1/2$. Then, by Lemma \ref{inversa}, we have that 
\[
\|g^{-1}-1\|<1, 
\quad \mbox{ so } \quad  \|\E(g^{-1})-1\|<1 ,
\] 
implying 
that $\E(g^{-1})$ is invertible. So, by Proposition \ref{lema1}, $g\in\ze$.

If $h\in\ze$, by the preceeding paragraph there is a neighborhood $V$ of $1$
such that $V\cap\ie\inc\ze$. Then $hV$ is a neighborhood of $h$ and, if 
$g\in hV$, then $h^{-1}g\in V$, so that $h^{-1}g\in\ze$. As $h\in\ze$, we have that 
$g\in\ze$, so $\ze$ is open. Clearly $g\ze$ is open for every $g\in\ie$, and as 
we can obtain $\ie$ as a disjoint union of sets $g\ze$, all open, they are
also closed.
$\ze$ can be easily seen to be connected, since it is the product of two 
connected groups. The last assertion becomes now clear.
\qed

\medskip

\noindent
We deduce that the group $\ze$ is closed, open and invariant in $\ie$. 
So we have a new Weyl group defined by 
\begin{equation}\label{w1e}
\wwe=\pi_0 (\ie) \simeq \ie / \ze
\end{equation}

\noindent
Next we shall see that the new Weyl group agrees with the old (unitary)
one. We first need the following Lemma: 
\smallskip

\begin{lem}\label{hene}
With the above notations, 
$$
\he= \ume \un = \nE\cap \gme \gn = \nE\cap\ze.
$$
\end{lem}
\dem Let $w\in\nE\cap\ze$. Then $w=mn$, with $m\in\gme,n\in\gn$. Using the 
polar decomposition of $m$ and $n$ and the fact that $\me \inc \n ' \cap \m$, 
we have
\[
w=v_m |m|\cdot v_n|n|=v_mv_n |m||n|=v_mv_n |mn|=v_nv_m |w|=v_nv_m,
\]
where $v_m\in\ume, v_n\in\un$, so that $w\in\he$.\qed

\medskip
\noindent
Now we have the technical tools we need to prove that both Weyl groups are
the same.

\begin{teo}\label{W(E)}
Let $\m$ be a \avn , $\E\in\em$. Then the Weyl 
group obtained by the unitary construction and the Weyl group 
obtained by the invertible construction are isomorphic (i.e. $\wwe=\we$).
\end{teo}
\dem Let $\varphi:\we\ra\wwe$ given by $\varphi([u]_{\we})=[u]_{\wwe}$, for 
$u \in \nE$. Then 
$\varphi$ is well defined and it is an isomorphism. 
Indeed, good definition is
clear since $\he\inc\ze$. By Lemma \ref{hene}, $\varphi$ is injective. 
Let $g\in\ie$. To see that the map $\varphi$ is onto,  
We must find a unitary $u\in\nE$ such that $[g]_{\wwe}=[u]_{\wwe}$.

\noindent
Since $g\in\ie$, we have that $(g^*)\inv \in\ie$ (by adjoining and using 
that $\E$ is $*$-linear) and so 
$g^*\in\ie$, since $\ie$ is a group. Therefore $g^*g\in\ie$ 
and, by Proposition \ref{lema1}, $g^*g\in\ze$. Then 
there exist $m\in\me$, $n\in\n$ with $g^*g=mn$. Using again the polar 
decompositions $m= v_m|m|$ and $n=v_n|n|$ with $v_m \in \ume$ and  
$v_n \in \un$, we have that 
\[
g^*g=mn= v_m|m|\cdot v_n |n|= v_mv_n |m||n|=
v_mv_n|mn|=v_mv_n|g^*g|=v_mv_n g^*g ,
\]
implying that $v_mv_n=1$, so we can write $g^*g=mn$ with $m\in\me^+$, 
$n\in\n^+$. Then $|g|=m^{1/2}n^{1/2}\in\ze$. Using the polar 
decomposition of $g$, there is a unitary $u\in\um$ with $g=u|g|$. As 
$|g|\in\ze$ and $g\in\ie$, it follows that $u\in\um\cap\ie=\nE$, so
$[g]=[u]$. Finnaly, since both groups are discrete, the mentioned 
isomorphism is also a homeomorphism. \qed

\ 

\begin{rem}\rm
Almost all the construction made in this paper can be extended trivially to
C$^*$ algebras. Proposition \ref{conexa} is the point were problems appear
since the invertible group of a C$^*$ algebra need not to be connected.
\end{rem}

\section{Differential geometry of $\oe$.}
In this section we shall consider only von Neuman algebras with separable 
predual, in order to assure the existence of faithful normal states. 

Let $\n\subseteq\m$ be von Neumann algebras.  $\E \in \emn $ and 
$\oe = \{g\E(g\inv\cdot g) g\inv : g\in \gm \}$.
The differential geometry of the orbit $\o$ has been already 
studied by Larotonda and Recht in \cita{pucho}, where it is assumed that 
$ \n '\cap\m \inc \n$.
In this case they show that $\oe$ admits a differentiable structure and 
the map $\Pi_\E :\gm \to \oe$ defines a reductive structure on $\oe$. 

The aim of this section is to remove that hypothesis, and we shall show that
the orbit $\o$ can be always given a differentiable structure, 
and even a unique reductive  structure if $\n$ is finite. We will also 
show that the existence of reductive structures for all conditional 
expectations $\E\in \em$ forces the algebra $\m$ to be finite.

\subsection{Differentiable Structure.}

Next we state some definitions and 
three classical Banach-Lie group theory's theorems that will be used 
afterwards. As a general reference about this subjet, see, for example, 
\cita {pucho2} or \cita {lang}.

\begin{fed}\rm
Given a Lie-Banach group $G$ (complex analytic, real analytic, or $C^\infty$), 
we denote by $L(G)$ the Lie algebra of $G$, 
which will be always identified (as a complex or real Banach space) 
with the tangent space  $T_1(G)$ of $G$ at the identity.
A subgroup $H$ of $G$ is called a {\bf regular} subgroup if it 
is also a Lie-Banach group (of the same type) and if $T_1H$ is closed and 
complemented in $T_1G$.

\end{fed}

\begin{teo}\label{lie2}
Let $G$ a Lie group, $H\inc G$ a subgroup such that there exist open sets $U,V$
with $0\in U$, $1\in V$ and a decomposition 
$T_1(G)=X\oplus Y$ (as a Banach space) 
satisfying
\ben
\item $\exp:U\ra V$ is a diffeomorphism
\item $H\cap V=\exp(X\cap U)$.
\een
Then $H$ is a regular subgroup of $G$ and $T_1(H)=X$.
\end{teo}

\begin{teo}\label{lie3}
Let $G$ be a Lie group, $H\inc G$ a regular subgroup. Then 
\ben
\item $G/H$ has a unique structure of differentiable manifold such that
$G\ra G/H$ is a submersion
\item $G\ra G/H$ is a principal bundle with structure group $H$
\item The action $G\times G/H\ra G/H$ is smooth.
\een
\end{teo}

\begin{teo}\label{lie4}
If $H$ is a subgroup of a Lie group $G$ and the connected component $H_1$ of
$1$ in $H$ is a regular subgroup of $G$, then 
$H$ is a regular subgroup of $G$ if and only if $H_1$ is open in $H$
\end{teo}

\vglue.4truecm
In the following Proposition  we construct a conditional expectation 
that will be essential in order to characterize the tangent space of $\oe$
(see also \cita {combes}).

\begin{pro}\label{ff} Let $\n\subseteq\m$ be von Neumann algebras 
and $\E\in \emn$. Fix a 
faithful normal state $\varphi$ on $\n$, and call $\psi=\varphi\circ\E$.
Then there exists a unique conditional expectation $\F\in {\cal E}(\m ,\me)$ 
such that $\E\F= \F\E$ and $\psi\circ\F=\psi$.
\end{pro}
\dem
Denote by $\sigma_t^\psi $, $t \in \zR$, the modular group of $\m$ induced by 
$\psi$. Since $\psi=\varphi\circ\E =\psi\circ\E  $, we have that 
$\sigma_t^\psi \circ \E = \E \circ \sigma_t^\psi$ for all
$t\in \zR$ (see \cita{combes} or \cita{take}). 
By direct computations we can deduce that
$\sigma_t^\psi(\me)=\me$ for every $t\in\zR$. Take $\F\in {\cal E}(\m ,\me)$ 
to be the unique expectation with $\psi\circ\F=\psi$
obtained by Takesaki's Theorem on the existence of conditional expectations 
\cita{take}. 
Since $\E |_{\me} \in {\cal E}(\me , \zn )$, then 
$\E\circ\F \in {\cal E}( \m , \zn)$ and $\psi\circ(\E\circ\F)=\psi$. 
When we represent $\m$ as usual in $L^2(\m,\psi)$, the three conditional 
expectations $\E,\F,\E\circ\F$ give rise to three orthogonal projections
$\e,\f,\g$ with $g=ef$. Since $g=g^*$, we have that $ef=fe$, so 
$\E\F=\F\E$.\qed

\vglue.4truecm
Using the expectation $\F:\m\ra\me$ from Proposition \ref{ff}, we can define
\begin{equation}\label{delta}
\Delta=E+F-EF \in {\cal B} (\m).
\end{equation}
$\Delta$ is a projection, since $\E$ and $\F$ commute. Its 
image is the closed subspace $\me+\n$ of $\m$, which 
can also be writen as a direct sum:
\[\Imagen\Delta=(\me\cap\ker\E)\oplus\n.\]

\begin{pro}\label{mdecomp}
With the preceeding notations, 
$\ze$ is a regular subgroup of $\gm$ and $T_1\ze
=(\me\cap\ker\E)\oplus\n$.
\end{pro}
\dem In order to use Theorem \ref{lie2}, we need a decomposition
\[T_1\gm=X\oplus Y\]
with $X=\me+\n$, the natural candidate to be $T_1\ze$. 
This decomposition exists because as $\gm$ is open in $\m$,
we have that $T_1\gm=\m$ and so the projection $\Delta$ introduced in the 
preceeding discussion  gives the 
desired decomposition. 

Note that the exponential map of the Banach-Lie group $\gm$ coincides with
the usual exponential map ($m \mapsto e^m$) under the identification of 
$L(\gm )$ with $\m$.  
As $\exp$ is a local diffeomorphism we can fix an open set $0\in U$ 
such that $\exp:U\ra V = \exp(U)$ is a diffeomorphism. Let $0 \in U' \inc U$ 
be an open set and $x\in U'\cap X$; 
then $x=a+b$ with $a\in\me$, $b\in\n$, so (as $a$ and $b$ commute), we 
have $\exp(a+b)=\exp(a)\exp(b)$ with $\exp(a)\in\gme$ and $\exp(b)\in\gn$, 
thus showing that $\exp(U'\cap X)\inc \exp(U')\cap\ze$. 

%If $y\in V\cap\ze$, then 
Let $0<\delta<1/2$ such that  
$$B_\m (1,\delta)=\{y\in \m : \|y-1\|< \delta\}\inc V. $$
Let $y \in B_\m (1,\delta)\cap\ze$.
%To see that $\exp^{-1}(y)\in U$ we must show that $y$ can 
%be written as a product of two elements in $\gme$ and $\gn$ respectively 
%such that they are both near $1$, so that we can assure that their preimages
%by $\exp$ are near $0$ and then their sum is in $U$. 
Let $g\in\me, h\in\n$ with $y=gh$. Note that $\F(h)$ is in $\gzn$. 
Indeed, since $h\in\n$, we have
that $\F(h)\in\zn$. To see that $\F(h)$ is invertible note that 
\[
\|g\F(h)-1\|=\|\F(gh-1)\|\le\|gh-1\|<\delta<1.
\] 
Now write $y=gh=(g\F(h))(\F(h)^{-1}h)$. Then 
\[\|g\F(h)-1\|<\delta\] as before and then, by Lemma \ref{inversa}, 
\[\|F(h)^{-1}g^{-1}-1\|<2\delta.\]
Notice also that $\|gh-1\|<\delta<1$ implies $\|gh\|<2$. Now,
collecting estimates, 
\[\begin{array}{rcl}
\|\F(h)^{-1}h-1\|&=&\|\F(h)^{-1}g^{-1}gh-1\| \\ && \\
&\le&\|gh\| \|\F(h)^{-1}g^{-1}-1\|+\|gh-1\| \\ && \\
&<&4\delta+\delta=5\delta.
\end{array}\]
Let $\eps>0$ small enough in order that $B_\m(0,2\eps)\inc U$. Let 
$\delta$ small enough such that 
\[\exp^{-1}(B_\m(1,5\delta))\inc B_\m(0,\eps).\]
Call $V'=B_\m(1,\delta) \inc V$, 
$U'=\exp^{-1}(V')\inc U$. Let $y\in V'\cap \ze$. Then $\exp^{-1}(y)\in U'$ 
and, since $g\F(h)$ and $\F(h)^{-1}h$ are in $B_\m (1,5\delta)$, 
their preimages $a=\exp^{-1}(g\F(h))\in \me$ and 
$b=\exp^{-1}(\F(h)^{-1}h)\in \n $ verify that $a+b \in U\cap X$. 
But $\exp (a+b) = y = \exp(\exp^{-1}(y))$ 
and $\exp$ is injective in $U$. Then $a+b =\exp^{-1}(y) \in U'$ and 
$$
\exp (U' \cap X) = V' \cap \ze 
$$

\qed

\begin{cor}\label{ieregular}
Let $\m$ be a \avn \ and  $\E\in\em$. Then, with the preceeding notations, 
the isotropy group $\ie$ is a regular subgroup of $\gm$.
\end{cor}
\dem We already know that $\ze$ is a regular subgroup, that $\ze$  is the 
connected component of $1$ in $\ie$ and that it is open in $\ie$. 
So by Theorem \ref{lie4}, $\ie$ is a regular subgroup.\qed

\bigskip

\begin{teo}\label{estructura}
Let $\m$ be a \avn \ and $\E\in\em$ a faithful normal conditional expectation.
Then the similarity orbit $\o\simeq\gm/\ie$, considered with the quotient
topology of the norm topology of $\gm$, can be given a 
unique complex analytic manifold structure such that it is 
an homogeneous space 
(i.e. the map $\Pi_\E:\gm\ra\o$ is a principal bundle with group 
structure $\ie$ and $\Pi_\E:\gm\ra\o$ is a 
submersion).
\end{teo}
\dem Apply Corollary \ref{ieregular} and Theorem \ref{lie3}.\qed
\bigskip
\begin{rem}\rm
The analogue of Theorem \ref{estructura} is still true (but replacing
complex analytic by real analytic) for the unitary
orbit $\uE \simeq \um / \nE $ 
under the action of the real analytic Banach-Lie group $\um $. 
\end{rem}

\input weillinois2.tex

%% file: def.tex
\makeatother
\def\dem{{\it Proof.\ }\rm}
\def\direcciones{\small \begin{tabular}{cccc}
Mart\'\i n Argerami&&&Demetrio Stojanoff\\
Dpto. de Matem\'atica&&&Instituto Argentino de Matem\'atica\\
Fac. Cs. Exactas -- UNLP&&&Saavedra 15 3er piso\\
(1900) La Plata&&&(1083) Buenos Aires\\
Argentina&&&Argentina\\
martin@mate.unlp.edu.ar&&&demetrio@mate.dm.uba.ar
\end{tabular}}

\newtheorem{fed}{Definition}[section]

\newtheorem{teo}[fed]{Theorem}
\newtheorem{lem}[fed]{Lemma}
\newtheorem{cor}[fed]{Corolary}
\newtheorem{pro}[fed]{Proposition}
\newtheorem{rem}[fed]{Remark}
\newtheorem{exa}[fed]{Example}
\newtheorem{num}[fed]{}

\def\zC{\hbox{\rm C\hskip -5.4pt\vrule height 8.0pt width 0.4pt depth -0pt\hskip 4.5pt}}
\def\zR{\hbox{$\displaystyle I\hskip -3pt R$}}

\def\[{$$}
\def\]{$$}
\def\inc{\subseteq}
\def\qed{\hfill\hskip 0.1truecm\vbox{\hrule\hbox{\vrule\hskip .2truecm
\vbox{\vskip .2truecm}\vrule}
\hrule}\hskip 0.1truecm}
\def\ben{\begin{enumerate}}
\def\een{\end{enumerate}}

\def\hrarr#1#2{\mathrel{ \mathop{\hbox to .5in{\rightarrowfill}} 
\limits^{\scriptstyle#1}_{\scriptstyle#2}  }}
\def\hlarr#1#2{\mathrel{ \mathop{\hbox to .5in{\leftarrowfill}} 
\limits^{\scriptstyle#1}_{\scriptstyle#2}  }}
\def\vdarr#1#2{\llap{$\scriptstyle#1$}\left\downarrow\vcenter to .5in{}\right.
\rlap{$\scriptstyle#2$}}

\def\ddrarr#1#2{\llap{$\vcenter {\hbox{$\scriptstyle #1$}}$}
\searrow\displaystyle\rlap{$\vcenter {\hbox{$\scriptstyle #2$}}$}}
\def\ddlarr#1#2{\llap{$\vcenter {\hbox{$\scriptstyle #1$}}$}
\swarrow\displaystyle\rlap{$\vcenter {\hbox{$\scriptstyle #2$}}$}}
\def\m{M}
\def\n{N}

\def\f{f}
\def\F{F}
\def\g{g}
\def\E{E}
\def\e{e}
\def\oe{{\cal S}(\E)}
\def\o{{\cal S}(\E)}
\def\emn{{\cal E}(\m,\n)}
\def\em{{\cal E} (\m)}

\def\um{{\cal U}_{\m}}
\def\un{{\cal U}_{\n}}
\def\me{\m_\E}
\def\ume{{\cal U}_{\me}} 
\def\gme{G_{\me}}

\def\G{L(G)}
\def\Q{{\cal Q}}

\def\gm{G_\m}
\def\gn{G_\n}
\def\gzn{G_{\zn}}
\def\ze{{\cal Z}_\E}
\def\ke{{\cal K}_\E}
\def\ind{\mbox{{\rm Ind}}(\E)}

\def\uE{{\cal U}(\E)}

\def\zn{{\cal Z}(\n)}

\def\zm{{\cal Z}(\m)}

\def\h{{\cal H}}
\def\he{\h_\E}

\def\hep{\h_{\E_p}}

\def\nE{{\cal N}_\E}

\def\nep{{\cal N}_{\E_p}}

\def\ie{I_\E}
\def\iF{I_\F}
\def\noi{}
\def\ra{\rightarrow}
\def\we{W(\E)}

\def\wwe{W_1(\E)}

\def\Imagen{\mbox{Im}\,}
\def\ee{{\cal X}(\E)}
\def\eee{{\cal Y}(\E)}
\def\eep{{\cal X}(\E_p)}

\def\inv{^{-1}}

\def\avn{von Neumann algebra}
\def\avns{\avn s}

\def\eps{\varepsilon}
\def\dpi{(\mbox{T } \Pi_\E)_{_1}}
\def\ad{\mbox{Ad}}
\def\z{{\cal Z}}
\def\up{{\cal U}(p)}

\def\cita#1{[\ref{#1}]}

%% file: weillinois2.tex
\subsection{Reductive Structure.}

Now we will start considering the conditions that will allow us to find a 
reductive structure in $\oe$ and to characterize it.
First we recall the definition of Homogeneous Reductive Spaces (see also
\cita {mr}):

\begin{fed}\label{HRS}
A Homogeneous Reductive Space (HRS) is a differentiable manifold $\Q$ and a 
smooth transitive action of a Banach-Lie group $G$ on $\Q$, 
$L:G\times\Q\ra\Q$ with:
\ben
\item {\bf Homogeneous Structure:} for each $\rho\in\Q$ the map
\[
\begin{array}{rcl}\Pi_\rho:G&\ra&\Q \\ g&\mapsto&L_g\rho \end{array}
\]
is a principal bundle with structure group $I_\rho=\{g\in G:L_g\rho=\rho\}$ 
(called the isotropy group of $\rho$).
\item {\bf Reductive Structure:} for each $\rho\in\Q$ there exists a 
closed linear subspace $H_\rho $ of the Lie algebra $\G$ of $G$ such that
$\G=H_\rho\oplus L(I_\rho) $ 
which is invariant under the natural action of $I_\rho$ and such that the
distribution $\rho\mapsto H_\rho$ is smooth.
\een
\end{fed}

In order to give a HRS structure to the orbit $\oe$ under the action of 
$\gm$, we must find a decomposition 
\[
L(\gm)=L(\ie)\oplus\ke,
\] 
such that the ``horizontal'' space $\ke$ verifies 
\begin{equation}\label{equiv}
g (\ke ) g\inv =\ke \quad \mbox{ for all } \quad g\in \ie. 
\end{equation}
Recall that $L(\gm)$ can be 
identified with $\m$, because $\gm$ is open in $\m$. Also, $L(\ie)$ can 
be regarded as $T_1\ie$, and also as $T_1(\ie)_1$ (where $(\ie)_1$ is the 
connected component of $\ie$ at $1$). We have shown in Proposition 
\ref{conexa} that the connected component if $\ie$ at $1$ is $\ze$, so 
$T_1\ie=T_1\ze=L(\ze)$, and by Proposition \ref{mdecomp}, $L(\ze)=\me+\n$. 
Therefore we already know (see equation (\ref{delta})) 
that such a decomposition of $\m$ can be found. The problem which arises now is
that we need a complement of $\me + \n$ verifying 
the equivariance property described in equation (\ref{equiv}).

\begin{lem}\label{dobleestrella}
Let $B\inc A$ algebras, $P:A\ra B$ a linear projection and 
$g\in G_A$ such that 
$g(\ker P)g^{-1}\inc\ker P$ and $ gBg^{-1}=B.$ Then 
$P(gxg^{-1})=gP(x)g^{-1}$ for every $x\in A.$
\end{lem}
\dem Straightforward.\qed

\begin{lem}\label{35}
Let $\n\inc\m$ \avns, $\E\in\emn$ a faithful normal conditional expectation.
Suppose that there exists  a  faithful normal {\bf tracial} state $\varphi$ 
of  $\n$. Let $\psi = 
\varphi \circ \E $ and $\F \in {\cal E}(\m , \me )$ as 
in Proposition \ref{ff}. 
Then 
\ben
\item The expectation $F$ is unique in the sense that 
for any other faithful normal tracial state $\rho$ in $\n$, the 
expectation $\F_\rho \in {\cal E}(\m , \me )$ induced by $\rho\circ E $ 
verifies $F_\rho= F$.
\item 
$\ie \inc \iF = \{g\in \gm : gF ( \cdot )g \inv = F( g \cdot g\inv ) \}.$
\een

\end{lem}
\dem 
We shall first show the uniqueness of $F$. 
For every faithful normal tracial state $\rho$ of $\n$, the corresponding 
$F_\rho \in {\cal E}(\m ,\me)$ given by Proposition \ref{ff} 
verifies that $F_\rho |_\n \in {\cal E} (\n , \zn )$ is the center valued 
trace of $\n$, since $\rho \circ F_\rho |_\n = \rho$ (see for example 8.3.10 of 
\cita {kadison}). Then 
$$
\psi \circ F_\rho = \psi \circ E \circ F_\rho = \psi \circ F_\rho \circ E  
= \psi \circ F_\rho|_\n \circ E  = \psi \circ F|_\n \circ E = \psi \circ F .
$$
So $F_\rho = F$. Fix now $g\in\ie$. It is easy to see that $g\me g^{-1}=\me$. 
Taking the polar decomposition $g = |g^*|u $, we know by Proposition \ref{lema1} 
and the proof of Theorem \ref{W(E)} that $u \in \ie \cap \um =
\nE$ and $|g^*| \in \ze$. Let us first see that $u \in \iF$. Indeed the expectation
$F_u = L_u (F )= u F( u^* \cdot u ) u^*$ verifies that 
$F_u \in {\cal E}(\m , \me )$,  $F_u \circ E = 
E \circ F_u $ and $\varphi ( u \cdot u^* ) \circ E \circ F_u = 
\varphi ( u \cdot u^* ) \circ E $. That 
is, $F_u$ is the expectation which corresponds by Proposition \ref{ff} to the 
trace $\varphi ( u \cdot u^* )$ of $\n$. By item 1 we can deduce that $F_u = F$, so
$u \in \iF$. Therefore it suffices to show that $\ze \inc \iF$ and, since 
$\ze = \gme \gn$, to show that $\gn \inc \iF$. Let $g \in \gn$, 
$y\in\me$ and $x\in\ker\F$. Then, using that  $\psi \circ \F = \psi$, we get
\[\begin{array}{rcl}
\psi(\F(gxg^{-1})y)&=&\psi(\F(gxg^{-1}y))=\psi(gxg\inv y)\\
&&\\
&=& \varphi(\E(gxg^{-1}y))= 
\varphi(g\E(xg\inv y g ) g\inv)\\
&&\\
&=& \varphi(\E(xg\inv y g )) = \psi(\F(xg^{-1}yg))\\
&&\\
&=&\psi(\F(x)g^{-1}yg)=0.
\end{array}\]
Then $\F(gxg^{-1})=0$, since $\F(gxg^{-1}) \in \me$ and $\psi$ is faithful. 
Thus  $g( \ker\F )g\inv \inc \ker \F $. By Lemma 
\ref{dobleestrella}, we conclude that $\ie \inc \iF $, showing item 2 \qed

\bigskip
\begin{pro}\label{HRS1}
Let $\n \inc \m$ be \avns , $\E\in\emn$ a faithful normal conditional
expectation and assume that $\n$ is finite. Then the similarity orbit 
$\oe$ has a unique  HRS structure under 
the action of $\gm$.
\end{pro}
\dem 
To find a reductive structure, we need to construct a decomposition 
$L(\gm)=L(\ie)\oplus\ke$, where $\ke$ is invariant by 
inner conjugation of elements of $\ie$. Fix a faithful normal tracial state 
$\varphi$ of $\n$ and consider $F\in {\cal E}(\m,\me)$ 
induced by $\varphi$ as in Proposition \ref{ff} and Lemma \ref{35}. 
By the discussion preceeding Proposition \ref{mdecomp}, 
it is clear that the projection 
$\Delta = I-(I-E)(I-F)$ gives the desired decomposition, i.e. 
$\ke=\ker\Delta$.

Now it remains to show that $\ie$ leaves $\ke$ invariant, and that the 
distribution $L_g \E\mapsto  g\ke$ is smooth. 
The first assertion follows, since 
$\ke=\ker\F\cap\ker\E$ and, by Lemma \ref{35}, $\ie \inc \iF$.

To see that the distribution is smooth, note that the projection 
onto $\ke$ with kernel $\me + \n$ is $I-\Delta = D= (1-\E)(1-\F)$. 
By Lemma \ref{35}, the map 
$$
\eta:\oe \to {\cal B}(\m) \quad \mbox{ given by }\quad 
\eta (\Pi_\E (g))= L_g D = (1-L_g \E)(1-L_g \F),  \quad g\in \gm
$$ 
is well defined and gives the desidered decomposition for all
$\Pi_\E (g) \in \oe$. Consider the commutative diagram
\[
\begin{array}{rcl}
\gm&\hrarr{Ad}{}&Gl({\cal B}(\m)) \\
\vdarr{\Pi_\E}{}&&\vdarr{}{\Pi_D} \\
\oe&\hrarr{\eta}{}&{\cal B}(\m)
\end{array}
\]
where  $\Pi_D(\alpha)=\alpha\circ D\circ\alpha^{-1}$, 
$\alpha \in Gl({\cal B}(\m))$. 
As we know that $\Pi_\E$ has analytic local cross sections by Theorem 
\ref{estructura}, the map $\eta$ is clearly analytic.

The uniqueness follows from the fact that our selection of $\ke$ (actually
the expectation $\F$) does not 
depend on the tracial state $\varphi$. Indeed, it is easy to see that 
for every faithful normal tracial state $\rho$ of $\n$, the corresponding 
$F_\rho \in {\cal E}(\m ,\me)$ given by Proposition \ref{ff} 
verifies that $F_\rho |_\n \in {\cal E} (\n , \zn )$ is the center valued 
trace of $\n$, since $\rho \circ F_\rho |_\n = \rho$. Then $F_\rho = F$ \qed
\bigskip
\begin{rem}\rm
Let $\n \inc \m$ be \avns , $\E\in\emn$ a faithful normal conditional
expectation and assume that $\me \inc \n$ even though $\n$ was not 
necessarily  finite. Then the assertion of Theorem \ref{HRS1} holds with the
same proof. Indeed, in this case $\ze = \gn$  
and one does not need a tracial state of $\n$ since $\Delta = \E$. 
This fact was already shown in \cita{pucho} under the slightly  more 
restrictive hypothesis that $\n ' \cap \m \inc \n$.
\end{rem}

\bigskip 

\begin {num}\label{matriz}\rm 
Let $\m$ be an infinite von Neumann algebra. Then there 
exists a properly infinite projection $p\in\zm$ such that $p\m$ is 
properly infinite and $(1-p)\m$ is finite. Let $\tau$ be a faithful normal 
trace in $(1-p)\m$. Since $p$ is properly infinite, 
it can halved, i.e. there exists a projection $q\in\m$ such that 
$q\sim p-q\sim p$, where $\sim$ denotes the von Neumann equivalence of 
projections. Using this projection $q$, we can identify $p\m$ with 
$q\m q\otimes M_2(\zC)$. So we identify $\m$ with $(q\m q)\otimes M_2(\zC))
\oplus (1-p)\m$.

Let $\n$ be the subalgebra $(q\m q\otimes 1)\oplus (1-p)\zC$ of $\m$. 
Consider the expectation $\E\in \emn $ given by 
\[\E=(\mbox{id}\otimes\mbox{tr}_2)\oplus\tau.\]
In matrices this can be seen as
\[
\E\left(\left(
\begin{array}{ccc}
a&b&0\\c&d&0\\0&0&x
\end{array}
\right)\right)
=\left(
\begin{array}{ccc}
\frac{a+d}{2}&0&0\\0&\frac{a+d}{2}&0\\0&0&\tau(x)
\end{array}
\right). 
\]
Straightforward calculations show that 
\[\n'\cap\m={\cal Z}(q\m q)\otimes M_2(\zC)\oplus (1-p)\m \]
and that $$\me=\n'\cap\m.$$

If $\oe$ admits an Homogeneous Reductive Structure, then there exists a 
bounded linear projection $P :\m\ra\n+\me$ with $g( \ker P )g\inv = \ker P$ 
for all  $g \in \ie$. 
Since $\un\i\ie$, then   
\[P(uxu^*)=uP(x)u^*\mbox{ for every }u\in\un \]
by Lemma \ref{dobleestrella}. 
Note that, as $(\n+\me)^* = \n +\me$, 
$P$ can be assumed to be $*$-linear. Indeed, if $P$ is not 
$*$-linear, we can replace it by 
\[P'(x)=\frac{1}{2}(P(x)+P(x^*)^*), x\in\m.\]
This $P'$ is also a projection onto $\n+\me$ and 
\[P'(uxu^*)=uP'(x)u^*\quad \mbox{ for every } \quad u\in\un.\]
\bigskip
\noi Since we know that
\begin{equation}\label{nemasn}
\n+\me=\left\{ \left(\begin{array}{ccc}n&z_2&0\\z_3&n+z_1&0\\0&0&m\end{array}\right):
n\in q\m q, z_i\in {\cal Z}(q\m q), m\in(1-p)\m\right\},
\end{equation}
it is clear that the elements in coordinates 2 1 and 1 2 of the image of $P$
belong to ${\cal Z}(q\m q)$.
Consider the linear map $T:q\m q\ra{\cal Z}(q\m q)$ given by
\[T(n)=\frac{1}{2}\left(P\left(\begin{array}{ccc}0&n&0\\n&0&0\\0&0&0
\end{array}\right)_{21} + P\left(\begin{array}{ccc}0&n&0\\n&0&0\\0&0&0
\end{array}\right)_{12}\right),  \ n \in q\m q\]
where $(\cdot)_{21}$ and $(\cdot)_{12}$ mean the 2 1 and the 1 2 
coordinate of the matrix respectively. 
Now we show the properties of $T$ that we will be of interest to us: 
\end{num}
\begin{pro}\label{delaf} 
Let $\m$ be an infinite \avn . Let $p$, $q$, $\n $ and $\E\in\emn $ as 
in \ref{matriz}. Assume that the orbit $\oe$ admits an Homogeneous 
Reductive Structure. Consider the linear maps $P$ and 
$T$ as before. Then the following properties are satisfied:
\ben
\item $T:q\m q\ra{\cal Z}(q\m q)$ is a $*$-linear mapping;
\item $T$ is a projection onto  ${\cal Z}(q\m q)$;
\item if $u\in{\cal U}_{q\m q}$, then $T(unu^*)=T(n)$ for every $n\in q\m q$;
\item $T(xy)=T(yx)$ for every $x,y\in q\m q$.
\een
\end{pro}
\dem 
\ben
\item That the image of $T$ is in ${\cal Z}(q\m q)$ can be seen in 
equation (\ref{nemasn}). $*$-linearity 
is clear since we  assume $P$ to be $*$-linear.
\item If $s\in{\cal Z}(q\m q)$, then the matrix
\[ \left(\begin{array}{ccc}0&s&0\\s&0&0\\0&0&0\end{array}\right)\]
is clearly in $\n+\me$, so the matrix is left invariant by $P$, and 
$T(s)=s$.
\item\label{uu} Let $u \in {\cal U}_{q\m q}$ and consider 
\[U = \left(\begin{array}{ccc}
u&0&0\\0&u&0\\0&0&1\end{array}\right) \in \un \inc \ie .
\] 
The basic property of $P$ is that $P(UmU^*) = UP(m)U^*$ for every 
$m\in \m$. But this clearly implies that  $T(unu^*)=uT(n)u^*
=T(n)$ for every $n\in q\m q$.
\item Follows from $\ref{uu}$ since the unitaries generate the whole algebra.
\een \qed

\begin{teo}\label{Principal}
Let $\m$ be a von Neumann algebra. Then the following conditions are 
equivalent:
\ben
\item The similarity orbit 
$\oe$ of any expectation $\E\in\em$ can be given a HRS structure under 
the action of $\gm$.
\item $\m$ is a finite von Neumann algebra.
\een
\end{teo}
\dem Let $ p$ be the biggest projection in $\zm$ such that $p\m$ is properly 
infinite, and $q$ a subprojection of $p$ that halves $p$, 
that is $q\sim p-q\sim p$. We shall use the 
notations of \ref{matriz} and the conditional expectation
$\E\in \em$ considered there. If condition 1. holds, using 
Proposition \ref{delaf} and \ref{matriz}, we can construct a 
``tracial'' bounded projection $T:q\m q\ra {\cal Z}(q\m q)$. 
Since $q$ is also 
properly infinite,  there is a projection $r\in q\m q$ such that  
$r\sim q-r\sim q$ in $q\m q$. 
Using ``traciality'' of $T$,
\begin{equation}\label{pyf}
T(q)=T(q-r)=T(q)-T(r)=T(q) -  T(q)=0.
\end{equation}
Recall from proposition \ref{delaf} that $T(q)=q$, 
so by equation (\ref{pyf}) we have that $q=0$, and this implies that $p=0$. 
So $\m$ is a finite von Neumann algebra.

Conversely, suppose that $\m$ is finite and $\E\in \em$. Then 
$\n = \E(\m)$ is a finite \avn \  and we can apply Proposition
\ref{HRS1}. \qed

\begin{rem}\rm
Let $\n \inc\m $ \avns \ and $\E\in \emn$ such that $\o$ has a 
structure of HRS. We shall describe explicitely the geometrical invariants
of $\o$. First we compute the tangent map at $1$ 
of the fibration $\Pi_E : \gm  \to \o$. 
For simplicity we shall consider $\o \inc {\cal B}(\m) $, 
in spite of the fact that the topology considered in $\o$ is not in general 
that induced 
by ${\cal B}(\m) $. In this sense, for $ x \in \m $, 
$$
\dpi (x) = [x,\E(\cdot)]-\E([x,\cdot]) ,
$$
where $[x,y] = xy-yx$, for $x,y \in \m$. Indeed,  
let $x\in \m$ and consider the curve $\alpha(t)=e^{tx}$. Note that 
$\alpha(0)=1$ and $\dot{\alpha}(0)=x$. Then 
\[
\begin{array}{rcl}
\dpi (x) &=& \frac{d}{dt}(\Pi_\E(e^{tx}))|_{t=0}\\
&&\\
&=&\frac{d}{dt}\left(\ad (e^{tx})\circ
\E\circ\ad (e^{-tx}) \right)|_{t=0}\\ 
&&\\
&=&\left((\ad (e^{tx}))'\circ\E\circ\ad (e^{-tx})+\ad (e^{tx})\circ\E\circ
(\ad (e^{-tx}))'\right)|_{t=0}\\ 
&&\\
&=&\left((\ad (e^{tx}))([x,\E\circ\ad (e^{-tx})])+\ad (e^{tx})\circ\E\circ
(\ad (e^{-tx})([\cdot,x])\right)|_{t=0}\\ 
&&\\
&=&[x,\E(\cdot)]-\E([x,\cdot]).
\end{array}\]
An interesting computation using this formula shows, as it must 
be, that $$\ker \dpi = \me+\n = L(\ie ) .$$
\smallskip 
On the other hand, if $\ke = \ker \Delta $ 
is the horizontal space at $\E$ of $\o$, 
then 
$$\dpi |_{_{\ke}} : \ke \to T(\o )_\E 
$$ 
is an isomorphism. It is usual consider its inverse 
$K_E :  T(\o )_\E  \to \ke $ in order to identify tangent vectors with
elements of $\m$ (see, for instance, \cita{mr}). With this convention we shall 
describe the torsion and 
curvature tensors, $T$ and $R$, respectively.
Let $V, W $ and $Z\in T(\o )_\E$. Then
\smallskip
\ben
\item $T(V,W)=\dpi([K_\E(V), K_\E(W)])$.
\item $R(V,W)Z=\dpi([K_\E(Z),\Delta([K_\E(V),K_\E(W)]))$.
\item The unique geodesic $\gamma$ at $\E$ such that $\dot{\gamma}(0)=V$
is given by 
$$
\gamma(t)=L_{e^{t K_\E(V)}}\E .
$$
\item The exponential map of $\o$ is given by 
$$
\exp_\E(X)=L_{e^{K_\E(X)}}\E \quad  \mbox{ for } \quad X \in T(\o )_\E .
$$
\een

\end{rem}

\section*{ References}

\begin{enumerate}
\item\label{acs1} E. Andruchow, G. Corach and D. Stojanoff; A geometric
characterization of nuclearity and injectivity, J. Funct. Anal. {\bf 133}
(1995), 474-494. 
\item\label{acs3} E. Andruchow, G. Corach and D. Stojanoff, 
Geometry of the sphere of a Hilbert module,   
Math.  Proc. of the Cambridge Math. Soc., to appear.
\item\label{alrs}E. Andruchow, A. R. Larotonda, L. Recht and D. Stojanoff; 
Infinite dimensional homogeneous reductive spaces and finite index
conditional expectations, Illinois J. of Math. {\bf 41} (1997), 54-76.
\item\label{as}E. Andruchow, D. Stojanoff, 
``Geometry of conditional expectations and finite index'', 
Int. Journal of Math. {\bf 5} (1994), 169-178.
\item\label{av2} E. Andruchow and A. Varela, 
Homotopy of State orbits, preprint (1998).
\item\label{ass}  H.Araki, M.Smith and L. Smith, On the homotopical 
significance of the type of von Neumann
algebra factors, Commun. Math. Phys. {\bf 22} (1971),71-88.
\item\label{ars}M. Argerami, D. Stojanoff, ``The Weyl Group and the 
Normalizer of a Conditional Expectation'', Int. Eq. and Op. Theory, 
34 (1999), 165-186.
\item\label{ars2}M. Argerami, D. Stojanoff, ``Some examples of 
the Weyl Group  of Conditional Expectations'', preprint.
\item\label{combes}F. Combes, C. Delaroche, 
``Groupe modulaire d'une esperance 
conditionnelle dans une algebre de von Neumann'', 
Bull. Soc. Math. France {\bf 103} (1975), 385-426.
\item\label{connes}A. Connes, ``Une classification des facteurs de type III'', 
Ann. Ec. Norm. Sup. Paris {\bf 6} (1973), 133-252.
\item\label{cpr2} G. Corach, H. Porta and L. Recht, Differential geometry 
of systems of projections in
Banach algebras, Pacific J. Math. {\bf 143} (1990), no. 2, 209-228.
\item\label{cpr} G. Corach, H. Porta and L. Recht, The geometry of spaces of 
projections in C*-algebras, Adv. in Math. {\bf 101} (1993), 59-77.
\item\label{greenberg}M. J. Greenberg, ``Lectures on Algebraic Topology'',
W. A. Benjamin, inc., New York, 1967. 
\item\label{h}J.-F. Havet, ``Esperance conditionelle minimale'', 
J. Oper. Th. {\bf 24} (1990), 33-55.
\item\label{han}  D. E. Handelman, K$_0$ of von Neumann algebras and
AFC$^*$-algebras, Quart. J. Math. Oxford (2) {\bf 29} (1978), 429-441.
\item\label{jones}V.F.R. Jones, ``Index for subfactors'', 
Invent. Math., {\bf 72} (1983), 1-25.
\item\label{kadison}R. Kadison, John R. Ringrose, ``Fundamentals of 
the Theory of Operator Algebras'', II, Academic Press, New York 1986.
\item\label{kosaki2}H. Kosaki, ``Characterization of Crossed Product
(Properly Infinite Case)'', Pacific J. of Math., {\bf 137} (1989), 159-167.
\item\label{lang}S. Lang, {\it Differentiable Manifolds}, Addison--Wesley, 
Reading, Mass. 1972
\item\label{pucho2} A. R. Larotonda, {\it Notas sobre variedades diferenciales,
Notas de Geometr\'\i a y Topolog\'\i a 1},
INMABB-CONICET, Universidad Nacional del Sur, Bahia Blanca (1980)
\item\label{pucho}A. Larotonda, L. Recht, 
``The orbit of a conditional expectation as a reductive homogeneous space'',
Monteiro, Luiz (ed.), Homage to Dr. Rodolfo A. Ricabarra. Bahia Blanca: 
Univ. Nacional del Sur, 61-73 (1995).
\item\label{longo-1}R. Longo, ``Index of subfactors and statistics of 
quantum fields I and II'', Comm. Math. Phys. {\bf 126} (1989), 217-247, 
{\bf 130} (1990), 285-309.
\item\label{mr}   L. Mata Lorenzo, L. Recht, Infinite
dimensional homogeneous reductive spaces, Reporte 91-11, U.S.B., (1991).
\item\label{popa-1}M. Pimsner, S. Popa, ``Entropy and index for 
subfactors'', Ann. scient. Ec. Norm. Sup., 4ta serie, {\bf 19} 
(1986), 57-106.
\item\label{Po}S. Popa, ``Classification of subfactors and their
endomorphisms'', AMS series CBMS  {\bf 86} (1995).
\item\label{PR}H. Porta and L. Recht; Minimality of geodesics in Grassmann
manifolds, Proc. Amer. Math. Soc., 100, (1987), 464-466. 
\item\label{sch} H. Schr\"{o}der, On the homotopy type of the regular 
group of a W$^*$-algebra, Math. Ann. {\bf 267} (1984), 271-277.
\item\label{take} M. Takesaki, 
Conditional Expectations in von Neumann Algebras, 
J. Funct. Anal. {\bf 9} (1972) 306-321.

\end{enumerate}

\vglue1truecm

\direcciones
\bigskip

\noi
AMS Classification Numbers: 46L10 and 46L99.

\end{document}